\documentclass[reqno,a4paper,11pt]{amsart}
\usepackage{mathtools} 
\usepackage{amsmath,amsthm,amsbsy,amssymb, comment, booktabs, pstricks}
\usepackage{txfonts} 

\usepackage[hyphens]{url}

 \usepackage[backref=page]{hyperref}
\hypersetup{
    colorlinks = true,
linkcolor={red},
urlcolor={blue},
citecolor={green},    
urlcolor = {blue},
citebordercolor = {0.33 .58 0.33},
 linkbordercolor = {0.99 .28 0.23},
 breaklinks=true
}
\usepackage{cite}

\theoremstyle{plain}
\newtheorem{theorem}{Theorem}[section]
\newtheorem{corollary}{Corollary}[section]

\newtheorem{problem}{Problem}[section]
\newtheorem{question}{Question}[section]
\newtheorem{conjecture}{Conjecture}[section]

\newtheorem*{Problem2b}{Problem\,2b}
\newcommand{\problembb}[1]{\mbox{}}

\theoremstyle{definition}

\theoremstyle{remark}

\newcommand{\cP}{\mathcal{P}}
\newcommand{\cL}{\mathcal{L}}
\newcommand{\cA}{\mathcal{A}}
\newcommand{\cS}{\mathcal{S}}

\newcommand{\cM}{\mathcal{M}}

\newcommand{\dsp}{\hbox{\kern0.2em}}

\newcommand{\ZZ}{\mathbb{Z}}

\newcommand{\NN}{\mathbb{N}}

\usepackage[pdftex]{graphicx}
\graphicspath{{IMAGES/}{./}}
\usepackage[skip=2pt,font={small},textfont=it]{caption}
\usepackage{subfigure}  

\begin{document}
\title[DOI$^2$]{DOI$^2$}

\author[Cristian Cobeli]{Cristian Cobeli}
\address{"Simion Stoilow" Institute of Mathematics of the Romanian Academy, 21 Calea Grivitei 
Street, P. O. Box 1-764, Bucharest 014700, Romania}
\email{cristian.cobeli@imar.ro}

%
%


\dedicatory{This  manuscript is based on the author's presentation  at\\
The Ninth Congress of Romanian Mathematicians,\\
June 28 - July 3, 2019, Gala\c ti, Romania.
}

\begin{abstract}
We discuss some seemingly unrelated observations on integers, whose close or farther away neighbors
show a complex of combinatorial, ordering, arithmetical or probabilistic properties, emphasizing 
puzzlement in more common expectations.

\end{abstract}
\maketitle

{\small 
The unusual short title 'DOI$^2$' requires an explanation. It is a recursive acronym of the parity 
adjectives odd and even (the sound of pronunciation) or just the uppercase writing of \textit{doi}, 
the Romanian word for two. 
In our manuscript, $2$ is the main character, or the \textit{red string} of the $1$st of March 
celebrating symbol named in Romanian '\textit{m\u ar\c ti\c sor}' (little March), which is embraced 
by the \textit{white string} embodying the parity of '\textit{noted divisors}'. Their double bow  
twinkles with the hanging tassel spread around the $2\times 2$ themes we chose to present here.
}

\section{From an elementary Olympiad problem to a Waring problem question}\label{Section1}

As sometimes happens in a matter intended to avoid any involvement of chance, 
at one of the many local math competitions, which is of major importance for the young beginners 
and their supporters, at the first stage of The National Mathematical Olympiad, Prahova county 
Romania, February 24, 2019, 
the written problem looked different from what the author probably wanted to ask the participants.

\begin{Problem2b}[G. Achim, ONM Prahova 2019]
\problembb{\mbox{}} \label{labelProblem2b}
Find the non-negative distinct integers $m,n,p,q$ knowing that
{\small \normalfont [RO: Determina\c ti numerele naturale 
nenule $m,n,p,q$ distincte, \c stiind c\u a]}
\begin{equation}\label{eqW}
m^3+n^3+p^3+q^3=10^{2021}.
\end{equation}
\end{Problem2b}

Leaving aside the necessity on the solutions to have distinct components, 
condition that does not oversimplify any approach to the problem,
the straight meaning of the requirement is to find all solutions of~\eqref{eqW}.
This is a classic Waring problem for cubes, with a particularly large constant term $N=10^{2021}$.

Originally, in 1770, Edward Waring, in his \textit{Meditationes Arithmeticae} (see~\cite{AK2018} 
for a recent cover), 
summarized his experiments, and adventured to state that any natural number can be written as a sum 
of at most $g(k)$ positive integers that are $k$-th powers. Thus, any natural number can be 
expressed as a sum of $4$ squares, or $9$ cubes, or $19$ fourth powers, and so on, he wrote.
A text in a nutshell that developed into a lot of dreams and endeavors for centuries.
A major step towards understanding the problems raised by Waring was made by Hilbert, who proved 
in 1909 the existence of $g(k)$ for all $k$.
%
Precisely, he showed that for any $k\ge 2$, there is a positive integer $g=g(k)$ such 
that any integer $n\ge 0$ can be represented as  $n=a_1^k+\cdots +a_g^k$,
for some natural numbers $a_1,\dots,a_g$.
As an aside, we mention that one might find it interesting to know more about the 
peaking events of that time in the mathematical life of G\"ottingen, as described in the fascinating 
biography written by Reid~\cite[XIV Space, Time and Number]{Reid1996}, including the touching fate 
of Minkowski, in whose memory Hilbert dedicated his work on Waring's problem.
Furthermore, the interested reader might start consulting 
the survey of Vaughan and Wooley~\cite{VW2002} and 
the articles of Deshouillers et al.~\cite{DHL2000}, 
Pollack~\cite{Pollack2011} and Siksek~\cite{Sik2016} for a summary of the results and methods 
developed over the years.
Two examples from the latest results that are valid for all $n$ except just a few 
particular cases are the following. Siksek~\cite[Theorem 1]{Sik2016} showed that 
any $n\ge 455$ can be written as a sum of seven cubes, making effective 
the same statement proved by Linnik~\cite{Lin1943} to hold true for $n$ sufficiently large, and
Deshouillers et al.~\cite[Conjecture 1,2]{DHL2000} who concluded their findings by conjecturing 
that there are exactly $113\,936\,676$ positive integers that cannot be written as a sum of four 
nonnegative integral cubes and the largest of them is the title of their article.


Turning back to Problem~\hyperref[labelProblem2b]{2b} and knowing the limited time the children 
have had available during the competition,  
while looking at how big $N$ is, one is driven to pinpoint two catches. 
First, there should be a concise way to write any solution and second, it should not be too 
difficult to enumerate all the solutions.
But instead, only an epsilon part of this is possible, unless one accepts as an answer 
to the Gordian Knot problem a G\"odelian-type presentation of all solutions in less than 100 words, 
without even a single one being able to be extracted explicitly.

Finding a solution is not too difficult using the hint in Problem 2a of the mentioned olympiad, 
which implied immediately 
that $1^3+2^3+3^3+4^3=10^2$. Then, observing that $2021-2=2019$ is divisible by $3$, one finds that 
the equality 
\begin{equation*}
   \begin{split}
   10^{3\cdot673}\left(1^3+2^3+3^3+4^3\right) = 10^{2019}\cdot 10^2
   \end{split}
\end{equation*}
is equivalent to
\begin{equation}\label{eqParticular}
   \begin{split}
   \left(10^{673}\right)^3 + \left(2\cdot 10^{673}\right)^3 +
   \left(3\cdot 10^{673}\right)^3 + \left(4\cdot   10^{673}\right)^3 = 10^{2021}\,,
   \end{split}
\end{equation}
that is, $S=\left(10^{673}, 2\cdot 10^{673},
   3\cdot 10^{673}, 4\cdot   10^{673}\right)$ is a solution 
of equation~\eqref{eqW}.
But next, finding other solutions or proving that $S$ is unique is not an easy task at all. 
Actually the facts are as follows.

\smallskip

Let $\nu(n)$ denote the number of representations of $n\in\NN$  as a sum of four cubes.
The problem raised later, after Waring, asked to show that $\nu(n) > 0$ for $n>n_0$.
For the upper bound of $\nu(n)$,
using the large sieve and bounds of exponential sums in several variables, 
Hooley~\cite{Hooley1978} showed that 
\begin{equation*}
   \begin{split}
   \nu(n) =O\left(n^{11/18+\epsilon}\right),
   \end{split}
\end{equation*}
while the expectation is that the true order of magnitude is around $O\big(n^{1/3}\big)$. 
Further, for the lower margin, Hooley~\cite{Hooley1977} has
showed that $\nu(n)$ is not of order  $o\big(n^{1/3}\big(\log\log n\big)^4\big)$.
Moreover, common experiments in Waring problems prove that $\nu(n)$ is within the above margins 
even if $n$ is relatively small.
Therefore, no question of uniqueness of the solution of \eqref{eqW}, but the number of solutions is 
huge, around $10^{673}$. Therefore, how could someone face the endeavor to write so many solutions 
during the allocated time of $2$ hours (reduced from $3$ hours, as it was in the old days), 
by comparing with some worldwide margins:


\begin{list}{$\circ$}{} 
\setlength{\itemindent}{-2mm}
\item
The number of atoms in the observable universe  is  $\approx 10^{81}$;
\item
The number of atoms in a A4-paper is   $\ll 10^{23}$;
\item
(Assume all atoms are grouped in A4-papers.) The number of A4-papers in universe    $\ll 
10^{58}$.
\end{list}


\noindent
Therefore, one needs to write about $10^{673}/10^{58}=10^{615}$ solutions on each A4-paper, 
which is quite a number compared with the number of characters of the King James authorized 
Bible, which is $3116480<10^7$ .

In spite of these big numbers, let us see that the idea of deducing solution~\eqref{eqParticular} 
can be exploited to find a \textit{tower of solutions} of~\eqref{eqW}.
Let $a,b\ge 0$ with $a+b=673$  be integers and let $(m_a,n_a,p_a,q_a)\in\NN^4$ be a solution of  
\begin{equation}\label{eqroot}
   \begin{split}
  \sum m_a^3 =10^{2+3a}.
   \end{split}
\end{equation}
Then $\sum m_a^3\cdot 10^{3b} =10^{2+3a+3b}= 10^{2021}$ and therefore,
$\big(m_a10^{3b}, n_a10^{3b},p_a10^{3b},q_a10^{3b}\big)$ is a solution of 
$\sum m^3 =10^{2021}$. This means that finding \textit{root solutions} $(m_a,n_a,p_a,q_a)$ of 
\eqref{eqroot}, equations in which the constant term is smaller, allows to find solutions of 
equation~\eqref{eqW}, where the constant term is large.
The smallest root solutions whose components are ordered lexicographically are the following: 

For $a=0,\ b=673$,  $1^3+2^3+3^3+4^3=10^2$.

For $a=1,\ b=672$, there are three  root solutions, for which: 
\begin{equation*}
   \begin{split}
6^3+24^3+34^3+36^3 =& 10^5,\\
10^3+20^3+30^3+40^3 =& 10^5,\\
12^3+16^3+34^3+38^3 =& 10^5\,.\
   \end{split}
\end{equation*}

For $a=2,\ b=671$ there are $43$ ordered root solutions listed in Table~\ref{TableRootSolutions}.

\begin{center}
\small
\begin{table}[h]
\caption{The lexicographically ordered solutions of the equation $x^3+y^3+z^3+v^3=10^8$.
}
\begin{tabular}{llll}
\toprule
(0,196,312,396)  & (44,64,250,438)  & (92,136,240,436)  & (155,309,322,322) \\
(4,122,295,417)  & (44,100,160,456) & (92,244,256,408)  & (156,176,244,424) \\
(4,302,304,354)  & (54,151,288,417) & (100,200,300,400) & (193,267,299,361) \\
(14,58,106,462)  & (58,134,256,432) & (100,256,272,396) & (200,210,295,385) 
\rule[1pt]{0mm}{5mm}\\
(18,107,220,445) & (58,159,337,386) & (107,184,213,436) & (204,256,292,368) \\
(18,200,232,430) & (58,188,319,393) & (114,147,277,420) & (216,260,298,358) \\
(22,263,316,369) & (60,240,340,360) & (114,170,274,418) & (225,295,300,330) \\
(28,44,358,378)  & (64,65,255,436)  & (120,160,340,380) & (230,288,295,337) 
\rule[1pt]{0mm}{5mm}\\
(32,124,148,456) & (67,92,352,381)  & (128,172,292,408) & (240,244,256,380)\\
(37,65,75,463)   & (70,183,198,441) & (145,170,340,375) & (260,265,274,351) \\
(41,57,79,463)   & (72,195,277,414) & (151,282,288,369) &                   \\
\bottomrule
\end{tabular}
\label{TableRootSolutions}
\end{table}
\end{center}

\noindent
Let us remark that the root solutions, whose components are ordered lexicographically, generate
$1\times 4!=24$ solutions (with no restrictions on the order of the non-negative integer 
components) of equation~\eqref{eqroot} with $a=0$ and
$3\times 4!=72$ solutions of equation~\eqref{eqroot} with $a=1$.
If $a=2$, excluding the $12$ permutations that would have been counted twice, since  
$(155,309,322,322)$ has tow equal components, we have 
$43\times 4!-12=1020$ solutions of equation~\eqref{eqroot} with $a=2$.
Also, compare with the expected order of magnitude in the last case, which is 
$\approx 10^{8/3}=464.15\dots$

Much harder is to find by brute force the next root solutions. For $a=3,\ b=670$, the first ordered 
solutions are:
\begin{equation*}
   \begin{split}
0^3+1960^3+3120^3+3960^3 =& 10^{11},\\
3^3+649^3+1775^3+4549^3 = & 10^{11}\,.
   \end{split}
\end{equation*}
Notice in the lists the towers of solutions being formed:
$(1,2,3,4)$; $(10,20,30,40)$;  $(100,200,300,400)$; and
$(6,24,34,36)$; $(60,240,340,360)$; and
$(0,196,312,396)$; $(0,1960,3120,3960)$.

As the power of the constant terms increases, it is more and more difficult to find corresponding 
new root solutions, despite their number increases fast. These solutions produce solutions of the 
original 
equation \eqref{eqW} with fewer and fewer  ending zeros. Thus, 
the unwished but inspiring 
formulation of Problem~\hyperref[labelProblem2b]{2b} raises the question 
of finding particular solutions of 
equation~\eqref{eqW}, since there are so many solutions, around $10^{673}$, but the probability to 
find one is very small, about $10^{673}/10^{2021}=10^{-1348}$.

\begin{problem}\label{Problem12}
Find a solution of equation $m^3+n^3+p^3+q^3=10^{2021}$ whose components do not end with no zeros. 
More generally, find solutions of general Waring problems whose components are not tower-wise 
related
to their constant term.
\end{problem}

\section{A covering problem and a question: Is it true that prime numbers are 
superabundant?}\label{Section2}

In this section we discuss a generalization of a problem of Pillai~\cite{Pillai1940} regarding 
finite sequences of consecutive integers with no member relatively prime to all the others. The 
problem appears in a question about the representation of a product of consecutive integers as a 
perfect power. 
Pillai treats the problem in a series of papers~\cite{Pillai1940, Pillai1941} and
Brauer~\cite{Brauer1940} is the first to prove completely the main conjecture of Pillai.
Later, Eggeleton~\cite{Eggleton1987} revisits Pillai's problem introducing a new language in terms 
of graphs and gives new proofs to the results of Pillai and Brauer. More recently,
Ghorpade and Ram\cite{GR2012} generalize the results for arithmetical progressions in integral 
domains. The theme was also met in several other contexts, whose links can be found by starting 
with the references within the cited papers.

\begin{figure}[ht]
 \centering
 \mbox{
 \subfigure{
    \includegraphics[angle=-90,width=0.40\textwidth]{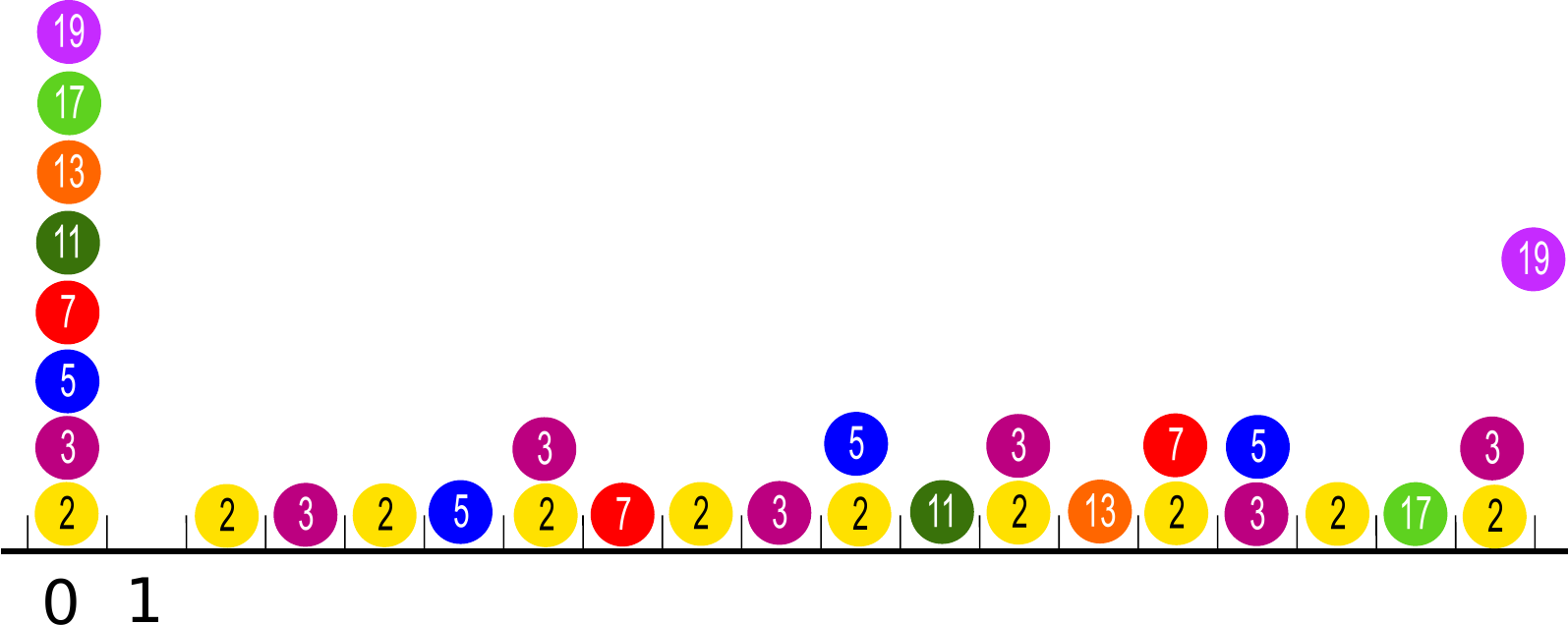}
 \label{CoveringAll}
 }\quad
 \subfigure{
    \includegraphics[angle=-90,width=0.45\textwidth]{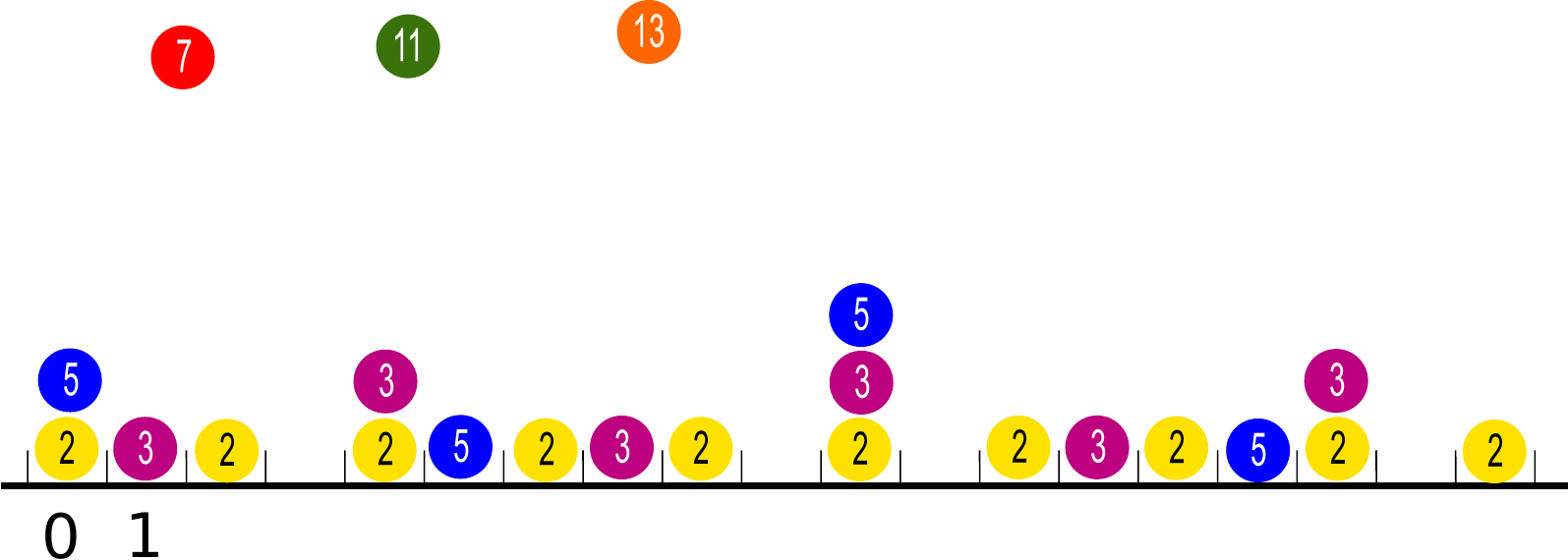}
 \label{Covering172}
 }
 }
\caption{
Two ways to fill the boxes of $\cL$ by pearls of primes. On the left, all primes are placed 
with a pearl in the zero box, exactly as in the sequence of integers, where  the box one always 
remains unoccupied,
while on the covering on the right side, primes $2$ and $5$ are 
also placed in the zero box, $3$ is shifted to the right and three more boxes are ready to fit 
effectively the remaining primes $7,11,13$ to generate the perfect minimal complete covering of the 
segment of length $17$. 
} 
 \label{FigCovering}
 \end{figure}

Let us consider the Eratosthenes sieve as a reverse process.  
Most of the discussion here can take place in a broader context, with numbers being not necessarily 
prime, but for brevity and simplicity we stick to the primes case.
Thus, let us say we have a large basket $\cP$ containing all prime numbers and an endless sequence 
of equal boxes situated on a straight line $\cL$, the analogue of the set of all integers. 
Further, let us think of any prime number as a sequence of pearls arranged on a wire at distances 
equal to its size [RO: \textit{\^in\c sir\u a-te m\u arg\u arite}]. 
The boxes of $\cL$ are supposed to be tall enough to fit any tower of pearls resulting from 
successive placement of primes on $\cL$.
During a placement of a prime $p$, when the left-right positioning of a prime is considered to be 
the suitable, the wire is dematerialized and the pearls of $p$ are attracted by gravity to the 
bottom.
In Figure~\ref{FigCovering} one can see two such ongoing arrangements.
 
We call \textit{segment} any finite sequence of consecutive boxes of $\cL$. The 
length of a segment $\cS$, which we denote by $l(\cS)$, is equal to the number of boxes of $\cS$.

Our interest will be to place \textit{effectively} a finite number of primes   on $\cL$ so that all 
boxes on a certain segment $\cS$ are filled. The \textit{requirement of effectiveness} demands that 
the pearls of any prime that is used in such an arrangement occupy simultaneously at least two boxes 
of $\cS$. This means that either $p$ is short enough, roughly less than half of $l(\cS)$, or $p$ is 
placed sufficiently close to the end points of $\cS$, so that two boxes of $\cS$ are filled by the 
pearls of $p$. 

We say that an arrangement of a set of primes $\cM\subset\cP$ \textit{covers} a segment 
$\cS\subset\cL$ if in each box of $\cS$ there is a pearl of a prime in $\cM$.
More precisely, in arithmetic language, identifying a segment of boxes filled by pearls of primes 
with a finite sequence of consecutive integers, we say that the segment is \textit{covered 
completely} if no element of the sequence is relatively prime with the product of all the others.

The description above is a natural solitary PL-game in which the player places prime pearls in 
straightly aligned boxes. The objective of the player is to fill effectively with pearls of primes 
as many boxes as possible of a certain  segment. The player wins if he succeeds 
to cover completely any segment of his choice, and his performance is the better the longer is the 
length of the covered segment.

A natural question to ask is whether there actually exist and how long are such segments that 
can be covered completely. The short answer is that there are no segments of length $\le 16$ that 
can be completely covered and any segment of length $\ge 17$ can be completely covered 
(see~\cite{Pillai1941, Brauer1940}).

\smallskip

Let us remark that the arrangements of primes that generate the sequence of all integers
(call it the $\ZZ$-arrangement) is close, but always fails to cover completely a segment with and 
endpoint at $0$ and, more generally, relatively close to $0$. In the $\ZZ$-arrangement all primes 
are placed calibrated with a pearl in the zero box, forming an infinitely high tower. Although,  
the neighbor boxes $\pm 1$ remained forever unfilled. At the other endpoint, if one wishes to 
extend any such almost completely covered segment, even the scantiest possible boxes always contain 
a marble of a prime with another marble in the zero box. 

A successful arrangement of primes on a completely covered segment of length $17$ is described by a 
sequence of $17$ consecutive integers $a,a+1,a+2,\dots a+16$, where $a$ is a solution of the 
following system of congruences:
\begin{equation}\label{System17}
   \left\{\begin{aligned}
         a    &\equiv  0\phantom{-} \pmod {2,5,11} \\
     a    &\equiv  -1 \pmod 3 \\
          a    &\equiv  -2 \pmod 7 \\
    a    &\equiv -3 \pmod {13} .\\          
   \end{aligned}\right.
\end{equation}
By the Chinese Remainder Theorem, the system~\ref{System17} has infinitely many solutions  
$a=27830+30030k, \ k\in\ZZ$. Similarly, any arrangement of a segment $\cS$ is reproduced 
infinitively many times on $\ZZ$, at places equally spaced in an endless arithmetic progression. 

An exhaustive analysis of all possibilities shows that the shortest complete covering is a segment 
of length $17$, and the arrangement is unique, except for its mirror, whose primes are placed in 
reversed order. The reversed covering actually has an initial positive solution of the 
system of congruences analogue to~\eqref{System17} that is closer to zero. This solution gives the 
completely covered segment with boxes labeled: $2184, 2185,\dots, 2200$, and the
periodicity of 
the solutions of the reversed covering is the same: $30030=2\cdot 3 \cdot 5\cdot 7\cdot 11\cdot 13$.

 \begin{figure}
{\footnotesize 
\begin{verbatim}
   b = -3  n = 27827 = 27827                  b = 10  n = 27840 = 2^6 * 3 * 5 * 29
   b = -2  n = 27828 = 2^2 * 3^2 * 773        b = 11  n = 27841 = 11 * 2531
   b = -1  n = 27829 = 17 * 1637              b = 12  n = 27842 = 2 * 13921
   -------------------------------------      b = 13  n = 27843 = 3 * 9281
   b =  0  n = 27830 = 2 * 5 * 11^2 * 23      b = 14  n = 27844 = 2^2 * 6961
   b =  1  n = 27831 = 3 * 9277               b = 15  n = 27845 = 5 * 5569
   b =  2  n = 27832 = 2^3 * 7^2 * 71         b = 16  n = 27846 = 2 * 3^2 * 7 * 13 * 17
   b =  3  n = 27833 = 13 * 2141              -----------------------------------------
   b =  4  n = 27834 = 2 * 3 * 4639           b = 17  n = 27847 = 27847
   b =  5  n = 27835 = 5 * 19 * 293           b = 18  n = 27848 = 2^3 * 59^2
   b =  6  n = 27836 = 2^2 * 6959             b = 19  n = 27849 = 3 * 9283
   b =  7  n = 27837 = 3^3 * 1031             b = 20  n = 27850 = 2 * 5^2 * 557
   b =  8  n = 27838 = 2 * 31 * 449           b = 21  n = 27851 = 27851
   b =  9  n = 27839 = 7 * 41 * 97            b = 22  n = 27852 = 2^2 * 3 * 11 * 211
\end{verbatim}
} 
\caption{The arrangement of prime pearls in the $17$ long completely covered segment. The 
factorization of the integers in the boxes b situated inside and near the minimal solution of the 
covering.}
\label{FigureFactorization17}
 \end{figure}


The fact that the complete covering~\ref{System17} of a segment of length 
$17$ can be extended naturally to longer segments, both to the left and to the right. This can be 
easily seen immediately in the list of the prime factorizations of the integers in the initial 
solution. In Figure~\ref{FigureFactorization17} one sees the barriers at boxes $\mathrm{b}=0$ and 
$\mathrm{b}=17$, which are occupied by large non-effective primes. But, we have the prime $17$ 
available, which can be placed in any of the two boxes to obtain a complete covering of length 
$18$. Furthermore, if for example, we place $17$ in box labeled $\mathrm{b}=17$, we are offered for 
free complete coverings of segments of length $18, 19, 20, 21$. Then, the new barrier can be 
overcome by placing $p=19$, which appears perfectly fit at our disposal. And the process can be 
continued finding sufficiently many available primes. 

For example, continuing this \textit{first come first served} arrangement process, one fills box
$\mathrm{b}=999$ by prime $p=647$ and has $50$ available primes remaining, and at a farther place
one arrives to fill box $\mathrm{b}=9191$ by prime $p=6043$ and remains with $351$ spare ones.

The question is: for how long this extension can be continued to the right by the first 
come first served rule of placing the smallest available prime in the new empty barrier box?

\begin{conjecture}\label{Conjecture17R}
Starting with the minimal covering of length $17$ and applying the 'first come first served rule', 
in which the first free box is filled by the smallest available prime, the extension extension of 
completely covered segments can be done indefinitely. Also, a similar fact should occur if one 
follows other more or less regular rules and by starting with different completely or almost 
completely covered segments.
\end{conjecture}
The 'first come first served rule' can be applied in different ways: always to the right of the 
segment to be extended, or always to its left. Also, other choices of the first empty 
box to fill, such as selecting alternatively from the right and from the left. Some other more 
complex rules of selection of the first box to fill seem, experimentally, to be even more 
efficient. We remark that in the process, the longer the segment, one may find several distinct 
complete coverings for~it. Almost always, these coverings may be taken as the starting point 
of Conjecture~\ref{Conjecture17R}, with still the same endless extension effect.
This phenomenon is similar to the generation of the Euclid-Mullin sequences
(see~Mullin~\cite{Mullin1963} and \cite{CZ2013}), but it offers easier more 
tangible results.

Furthermore, when the segment becomes long enough, we find seemingly that we can even skip some of 
the available odd primes by keeping them out for good of the arrangement process and without being 
too much disturbed in it. We don't know to prove that an infinite extension is possible by the 
``first come first served' rule for any starting complete covering, but we can show that there is a 
pattern of germ-coverings  that produce completely covered segments of any size. Moreover, these 
patterns create complete coverings of segments even when keeping out of the processes any finite 
set of odd primes. The lengths of completely covered segments are all integers greater than 
a certain size depending only on the set of primes left aside.

\begin{theorem}[\cite{CC2019}]\label{TheoremPL1}
 Let $k\ge 1$  and let $\cM=\{q_1,\dots,q_k\}$ be a finite set of odd primes. Then there exists 
$n_\cM\in\NN$ 
such that for any $N\ge n_\cM$ there exists a sequence of $N$ consecutive integers such 
that no one of 
them is relatively prime to the product of all the others. Moreover, the greatest common divisor of 
any number in the sequence and the product of all the others is divisible by a prime that is 
different to any $q\in\cM$. 
\end{theorem}

Theorem~\ref{TheoremPL1} is effective, meaning that the bound $n_\cM$ can be calculated 
explicitly, and it contains the conjecture of Pillai as a particular case.
\begin{corollary}(Brauer\cite{Brauer1940}) \label{CorollaryPL1}
 For any $N\ge 17$, there exists a sequence of $N$ consecutive integers such that no one of them is 
relatively prime to the product of all the others.
\end{corollary}
As an example of a choice of $\cM$ in Theorem~\ref{TheoremPL1},  if the the first odd prime, $q=3$, 
is excluded from the operation of filling the boxes, all segments of lengths greater than some 
integer smaller than $1300$ can be completely filled.

\begin{corollary}\label{CorollaryPL2}
 For any $N\ge 1300$, there exists a sequence of $N$ consecutive integers such that no one of them 
is relatively prime to the product of all the others, and moreover, the greatest common divisor of 
any 
integer in the sequence and the product of all the others divided by  the largest power of $3$ in 
its decomposition is $\ge 1$.
\end{corollary}

Using a different initial arrangement and taking account of the larger number of boxes that 
might not be filled by the pearls of the even prime $q=2$, one can still show that 
Theorem~\ref{TheoremPL1} holds true even with the restriction imposing the primes to be odd is 
removed.

\begin{theorem}[\cite{CC2019}]\label{TheoremPL2}
Let $k\ge 1$  and let  $\cM=\{q_1,\dots,q_k\}$ be a set of prime numbers. Then there exists 
$n_\cM\in\NN$ such that for 
any $N\ge n_\cM$ there exists a sequence of $N$ consecutive integers such that none of them is 
relatively prime to the product of all the others. Moreover, the greatest common divisor of any 
number in the sequence and the product of all the others is divisible by a prime that is different 
to any $q\in \cM$. 
\end{theorem}

From Theorem~\ref{TheoremPL2} it follows that the chances of a player to win remain intact, no 
matter his initial choices for a finite number of steps. The only possible inconvenience might be 
just the extension of the length of the segment that would ensure his win.

\begin{corollary}\label{CorollaryPL3}
 A player of the  solitary $\mathrm{PL}$-game has a strategy to win no matter what choices  he
followed for a finite number of moves.
\end{corollary}

\bigskip

\section{Constant digit numbers in the sequence of self-powers}\label{Section3}
Usually, according to the law of large numbers,
one expects that two sparse sequences of integers have few points of intersections.
This happens even when there is some regularity in the definition of the two sequences, but their 
nature is different.
Is it possible that the number of points of incidence is infinite?
In the following the two sequences we will consider are the sequence of self-powers and the 
sequence of integers whose representation in base $b$ (we restrict to the case $b=10$) has all 
digits equal.

For any non-negative integer $n$, let $l(n)$ denote the number of digits of $n$.
We say that $n\in\NN$  is a \textit{constant word number} if $l(n)$  has all digits equal.
Our startling example was noticed two years ago~\cite{CC2017,CZ2017} with the following 
match:
\begin{equation}\label{eq2017}
    l\left(2017^{2017}\right)=6666,
\end{equation}
so $2017^{2017}$ is a constant word number in base $10$. Remark that one needs some space to write 
the digits of $2017^{2017}$, since a dense A4 sheet may hardly contain $6000$ characters.

All small numbers are constant word numbers, since $l(n)$ has just one digit for all 
$n\le 999\hbox{\kern0.2em}999\hbox{\kern0.2em}999$. 
The next constant word number is $10^{10}$, which is the first member of a group that ends with 
$10^{11}-1$. They all have $11$ digits. Constant word numbers appear in groups that are longer and 
longer but farther and farther apart, like in a sort of generalized geometric progression. These 
groups are composed of numbers that have $1, 2,\dots, 9$, $11$, $22,\dots, 99$, $111, 222,\dots, 999$, 
$1111, 2222,\dots$ digits, respectively.
 
 In view of example~\eqref{eq2017}, one might wonder if there are other special years, for which 
their self-power is a constant word number. And the answer is that there are.
The previous one occurred $300$ years before, since
$1717^{1717}$ has $5555$ digits.
The next close misses are 2312 and 2602, since 
$l\left(2312^{2312}\right)=7778$ and
$l\left(2602^{2602}\right)=8887$. 
For the next real special years, we have to wait till 2889 and 3173, for which
 $ l\left(2889^{2889}\right)=9999$ and
$ l\left(3173^{3173}\right)=11111$.

The list of the constant word numbers starts with 
$1, 2,\dots,9, 10, 35, 46, 51, 194, 234, 273, 349, 423$ (see Table~\ref{TableCWN} for the number of 
digits of their self powers).

\begin{center}
\small
\begin{table}[h]
\caption{The first $23$ self powers constant word numbers.
}
\begin{tabular}{c|c|cccc|cccccc|cccc}
\toprule
$n$ & $1$-$9$ & $10$ & $35$ & $46$ &   
 $51$ & $194$ & $234$ & $273$ & $349$ & $386$ & $423$ & 1411 & 1717 & 2017 & 2889\\
$ l\left(n^n\right)$ & $1$ & $11$ & $55$ & $77$ & 
 $88$ & $444$ & $555$ & $666$ & $888$ & $999$ & $1111$ & 4444 & 5555 & 6666 & 9999\\
\bottomrule
\end{tabular}
\label{TableCWN}
\end{table}
\end{center}
And there are more such self power constant word numbers, for example
$n=631\dsp296\dsp394$, for which $l\left(n^n\right)$ has $5\dsp555\dsp555\dsp555$ digits.

\begin{question}\label{Question31}
Given the following two sequences: $\cS_1$ of the constant word numbers, and $\cS_2$ of the number 
of digits of self powers, that is,
\begin{equation*}
	\begin{split}
	\cS_1 &: 1,2,\dots, 9, 11, 22, 33,\dots,99, 111, 222,\dots,999, 1111, 2222,\dots\\
	\cS_2&: l\left(n^n\right),\ \text{for } n\ge 1,
	\end{split}
\end{equation*}
  how many common points do they have?
\end{question}
Checking the gaps, we find that two type of gaps that increase exponentially combine to separate 
the elements of $\cS_1$, while the average gap between the elements of $\cS_2$ is 
asymptotically equal to $e\cdot n^{n+1}$, for $n\ge 1$.

\medskip

We can also consider pairs of constant word self powers. We say that $m$ and $n$ are   
\textit{amicable constant word self powers} if $l\left(m^n\right)$ and $l\left(n^m\right)$ 
are constant word numbers. For example, such amicable pairs are:
$(4,368)$ since $4^{368}$ has  $222$ digits and $368^4$ has  $11$ digits;
$(39,698)$ since $39^{698}$ has  $1111$ digits and $698^{39}$ has  $111$ digits; and
$(48,66)$ since $48^{66}$ has  $111$ digits and $66^{48}$ has  $88$ digits.

More generally, we look at any size analogue of amicable pairs. Thus, for any  $k\ge 1$, 
we say that a tuple of positive integers 
$(m_1,\dots,m_k)$ is a tuple of \textit{amicable constant word of self powers} if each of 
the numbers
$l\left(m_1^{m_2}\right)$, 
$l\left(m_2^{m_3}\right)$,\dots, $l\left(m_{k-1}^{m_k}\right)$, and
$l\left(m_k^{m_1}\right)$
are written with only one digit. 
If $k=1$ the amicable $k$-tuples coincide with self powers constant word numbers.
Here are some examples of amicable tuples:

\medskip

\begin{list}{$\circ$}{} 
\setlength{\itemindent}{10mm}
     \item
$(26,62,49)$:\ \ $26^{62}$ $\leftsquigarrow$ ${\color{black}88}$ digits; 
$62^{49}$ $\leftsquigarrow$ ${\color{black}88}$ digits; 
$49^{26}$ $\leftsquigarrow$ ${\color{black}44}$ digits
     \item
$(49,39,62)$:\ \ $49^{39}$ $\leftsquigarrow$ ${\color{black}66}$ digits; 
$39^{62}$ $\leftsquigarrow$ ${\color{black}99}$ digits; 
$62^{49}$ $\leftsquigarrow$ ${\color{black}88}$ digits.
\end{list}

\medskip

\begin{list}{$\circ$}{} 
\setlength{\itemindent}{-4mm}
     \item
$(26,31,22,49)$:\ \ 
$26^{31}$ $\leftsquigarrow$ ${\color{black}44}$ digits; 
$31^{22}$ $\leftsquigarrow$ ${\color{black}33}$ digits; 
$22^{49}$ $\leftsquigarrow$ ${\color{black}66}$ digits;
$49^{26}$ $\leftsquigarrow$ ${\color{black}44}$ digits
     \item
$(66,54,25,47)$:\ \ 
$66^{54}$ $\leftsquigarrow$ ${\color{black}99}$ digits; 
$54^{25}$ $\leftsquigarrow$ ${\color{black}44}$ digits; 
$25^{47}$ $\leftsquigarrow$ ${\color{black}66}$ digits;
$47^{66}$ $\leftsquigarrow$ ${\color{black}111}$ digits.
\end{list}

\medskip

Counting $k$ amicable tuples with components less than a given fixed margin, we found that the 
chances to find an amicable $k$ tuple decrease with $k$, although the total number of $k$ amicable 
tuples increases significantly with $k$. Then, again, the basic question is whether there are 
really an infinite number of such amicable tuples.

\begin{question}\label{Question32}
How many amicable tuples of constant word self powers exist?
\end{question}

\section{Contrasting Shapes of Sturmian words}\label{Section4}

A common preconception at that initial shallow contact with a problem of a 
probabilistic nature regarding the parity of the involved objects is the expectation of a 
fifty-fifty occurrence of odds and evens. We present a case in which there is a clear bias between 
odds and evens, in fact the reflection of a plainly wide spread phenomenon regarding the parity of 
the divisors of the terms of sequences belonging to some quite different classes~\cite{CZ2018}.

Our object here are the 'simplest' non-finally periodic binary sequences. These are known as 
Sturmian words and the simplicity condition for a binary word $w$  is the requirement 
that $p_w(n)=n+1$ for all $n\ge 1$.  Here $p_w(n)$ is the complexity function of $w$, which by 
definition counts the number of distinct sub-words of length $n$. There is a gap from Sturmian words 
to the class of ultimately periodic sequence, whose characteristic is the fact that their complexity 
function is bounded.

Since their introduction by Hedlund and Morse~\cite{MH1940}, Sturmian words were intensely studied 
by different authors within a broad area of interests (see~\cite{AS2003,B1986, BS2002, LR2007,L2011,LL2015,M2004,MM2006, M2000, OEIS,RV2017}).

Following the ideas presented in~\cite{CZ2017} and~\cite{CZ2019},
we wish to show the two contrasting faces of Sturmian words, the more regular fractal-type face and 
the the random looking one. The asymptotic bias slope of the parity of the divisors 
function~\cite{CZ2017} and~\cite{CZ2019} is a type of phenomenon seen also in other contexts, such 
as those discussed in the following works:
\cite{CVZ2003,
CZ2001,
CZ2006,
CZ2016,
RS94,
L2012,
CGZ2003,
CVZ2012,
CZ2002}.

\subsection{Fractal face of Sturmian words}
A classic example of Sturmian words is the Fibonacci word, generalized by 
Dumaine~\cite{Dumaine2009}, and Ram\'irez et al.~\cite{Ramirez2014a, Ramirez2014b}.
They are generated using the concatenation operation in their defining recursive formula.

A general convenient way to define Sturmian words is using the \textit{rotation function} with two 
parameters modulo one.
Denote $R_\theta(\varphi)=\varphi+\theta\pmod 1$, where
$\varphi\in[0,1)$ 
and  suppose $\theta$ is irrational. 
Then define the word $w=w_1w_2\cdots$, with letters $w_n=a$ if 
$R_\theta^{(n)}(\varphi)\in[0,\theta)$ and $w_n=b$, else. For example:
\begin{equation}\label{eqsqrt7pi}
   \begin{split}
    R_{\sqrt{7}/7}{(0.2)} & \text{ generates } 
w=abbababbabbababbababbabbababbabbab\dots\\
    R_{\pi/8}{(0.2)}\phantom{R} & \text{ generates } 
w=abbababbababbababbababbababbabbaba\dots
   \end{split}
\end{equation}
Notice that both words in \eqref{eqsqrt7pi} contain exactly three distinct sub-words of length $2$ 
and $aa$ is not a sub-word, verifying the particular case of the complexity condition $p_w(2)=3$.

A standard drawing rule named the \textit{odd-even drawing rule}, was used 
originally to draw Fibonacci fractals~\cite{Dumaine2009, Hoffman2018}.
Starting at the origin and looking up, the curve associated to the Sturmian word 
$w=w_1w_2w_3\dots$ are constructed as follows:
\begin{itemize}
\setlength{\itemindent}{10mm}
 \item[$\circ$]
 $w_n=a$: walk one step forward;
 \item[$\circ$] 
 $w_n=b$ and $n$ even:  walk one step forward and turn left $90^\circ$;
 \item[$\circ$]
 $w_n=b$ and $n$ odd:  walk one step forward and turn right $90^\circ$.
\end{itemize}

\begin{figure}[ht]
 \centering
 \mbox{
 \subfigure{
    \includegraphics[width=0.48\textwidth]{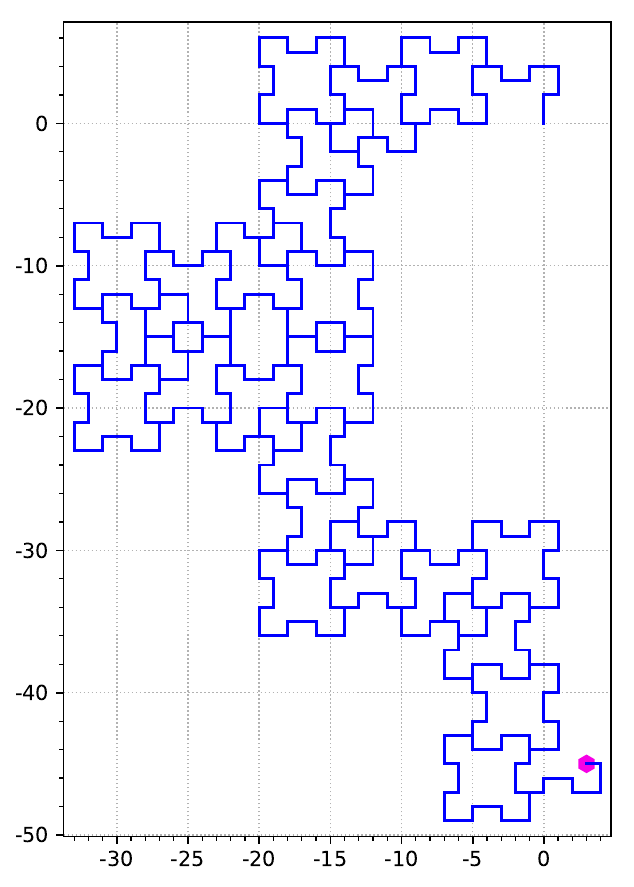}
 \label{Figsqrt7pe71000}
 }\quad
 \subfigure{
    \includegraphics[width=0.467\textwidth]{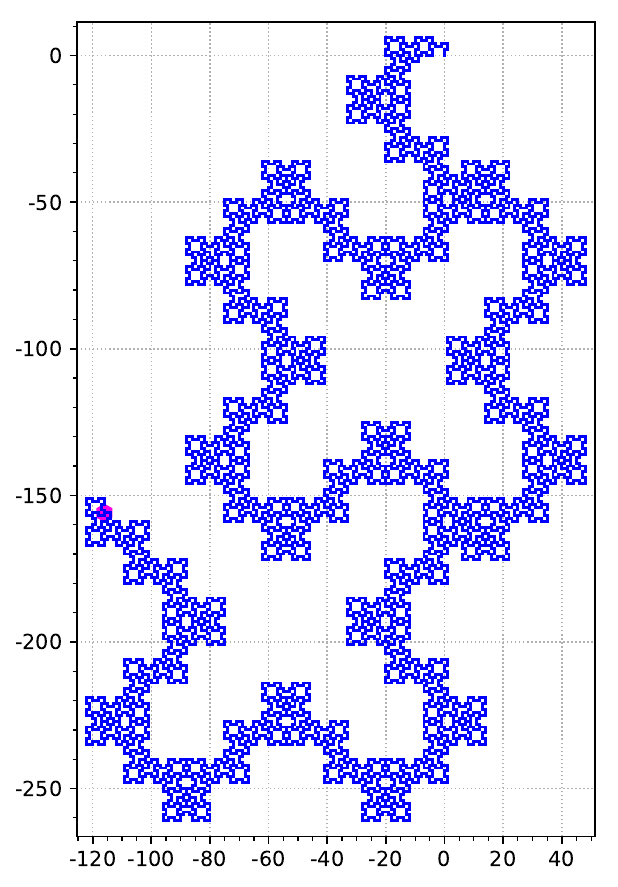}
 \label{Figsqrt7pe720000}
 }
 }
\caption{
Fractal type curves generated by the “odd-even drawing rule” applied to the Sturmian 
word defined by the rotation $R_{\sqrt{7}/7}{(0.2)}$. The image on the left represents 
the trajectory after $1000$ steps and the one on the right after $20000$ steps  (so the image on the 
left side can be found embedded in image on the right side).
%
} 
 \label{Figsqrt7pe7}
 \end{figure}

\begin{figure}[ht]
 \centering
 \mbox{
 \subfigure{
    \includegraphics[width=0.48\textwidth]{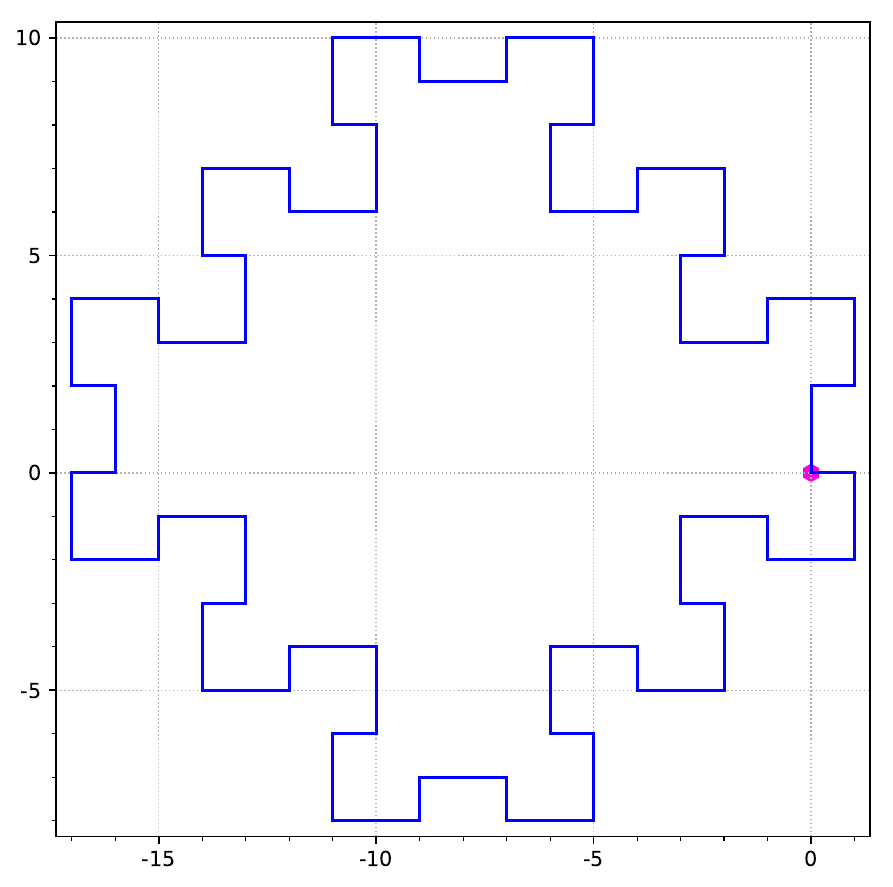}
 \label{Figpipe8112}
 }\quad
 \subfigure{
    \includegraphics[width=0.467\textwidth]{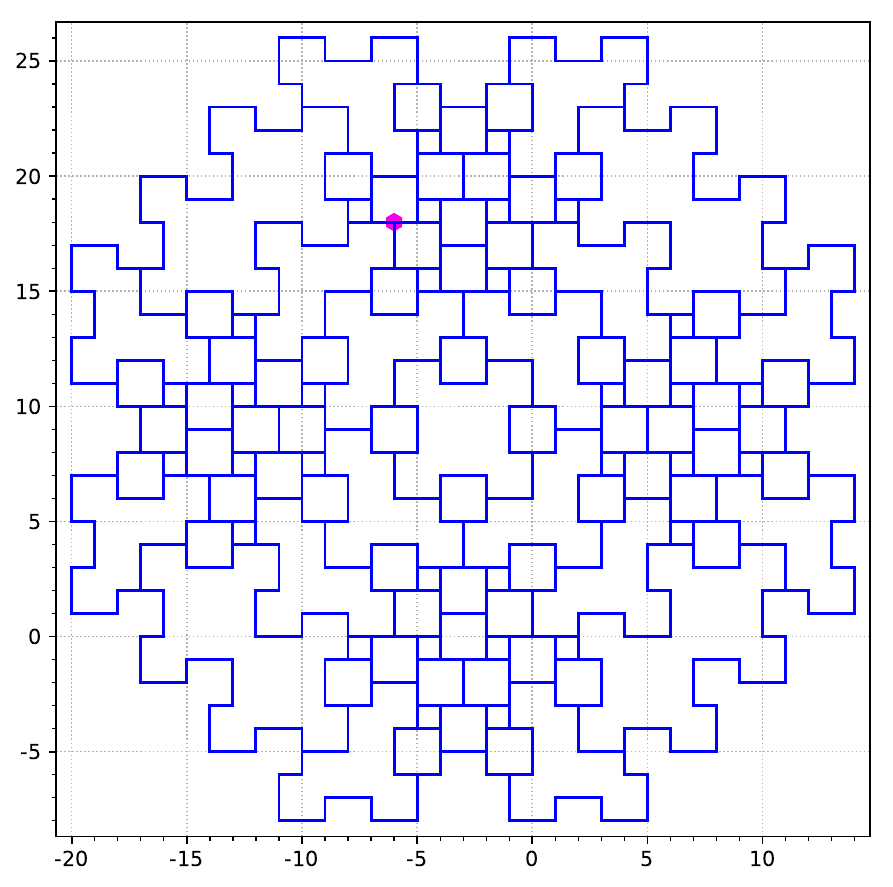}
 \label{Figpipe820000}
 }
 }
\caption{
Fractal curves generated by the “odd-even drawing rule” applied to the Sturmian 
word defined by the rotation $R_{\pi/8}{(0.2)}$. The image on the left side represents 
the trajectory after $200$ steps and the one on the right side is the continuation for a total of 
$20000$ steps.
%
} 
 \label{Figpipe8}
 \end{figure}

Applying the odd-even drawing rule to arbitrary Sturmian words we found a great variety of 
fractal-type shapes. Two examples are shown in Figures~\ref{Figsqrt7pe7} and ~\ref{Figpipe8}.
Let us mention that in various experiments we have noticed that there are many unexpected 
contrasting 
situations related to the complexity of the curves and that of the generators. For example, it is 
not rare to find simpler curves if the rotation parameter $\theta$ is transcendent than in the case 
where $\theta$ is algebraic.

 \subsection{Random walks of the divisor parity slope on Sturmian trajectories}
Fix the alphabet $\cA=\{a,b\}$ and let $w=w_1w_2\cdots$ be a word,  with letters $w_j\in\cA$.
We define the following counters of the parity of the divisors ranks:
\begin{equation*}
	\begin{split}
	o_w(n)&:=|\{j\in\NN : j \text{ divides } n, w_j=b \text{ and } n/j \text{ is odd } \}|,\\
	e_w(n)&:=|\{j\in\NN : j \text{ divides } n, w_j=b \text{ and } n/j \text{ is even } \}|\,.
	\end{split}
\end{equation*}
To illustrate the meaning of the above parity function, let $n\ge 1$, suppose $n=r\cdot s$ and make 
the association of the rank to the pair of the divisors: $n  \rightsquigarrow (r,s) $.
Then, the $r$--flag checks the letter $w_r$ (is it $a$  or $b$?)
and, if the flag is right, the parity of $s$ increments correspondingly either $o_w(n)$ or $e_w(n)$.
The calculation of $o_w(n)$ and $e_w(n)$ ends after all the pairs of the divisors of $n$ are 
evaluated.

\bigskip

\begin{figure}[ht]
 \centering
 \mbox{
 \subfigure{
    \includegraphics[width=0.40\textwidth]{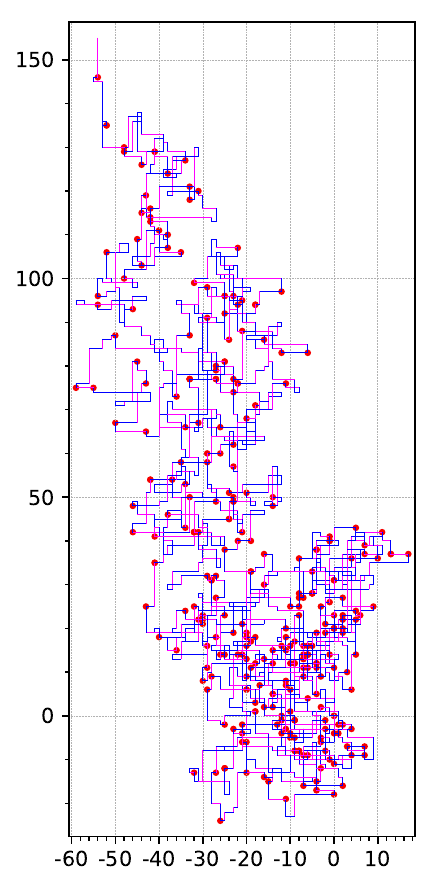}
 \label{FigRandomWalks77}
 }\quad
 \subfigure{
    \includegraphics[width=0.52\textwidth]{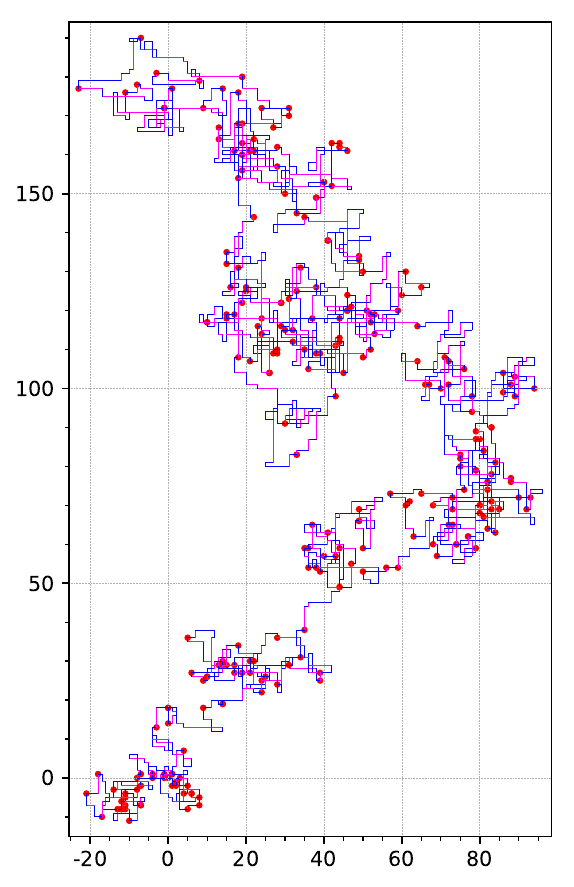}
 \label{FigRandomWalkspipe8}
 }
 }
\caption{
The first $2000$ steps of two random walks generated by the parity divisors slope of Sturmian 
words defined by rotations $R_\theta(\varphi)$, with $\varphi=0.2$ and $\theta=\sqrt{7}/7$ (left 
side) and  $\theta = \pi/8$~(right side). The red dots are drawn to indicate the zero-length steps 
(places where the odd and even divisor parity functions are equal).
} 
 \label{FigRandomWalks}
 \end{figure}

\textbf{Examples}: 
Let $w=abbaa\hspace{.5mm}bbbaa\hspace{.5mm}ba\dots$. Then, we have:

\medskip


\begin{minipage}[t]{.5\textwidth}
 $n=8 \rightsquigarrow (1,8)$ \quad $w_1=a$\ $\not\!\checkmark$ 

$n=8 \rightsquigarrow (2,4)$ \quad $w_2=b$\ $\checkmark$, $4$ even

$n=8 \rightsquigarrow (4,2)$ \quad $w_4=a$\ $\not\!\checkmark$ 

$n=8 \rightsquigarrow (8,1)$ \quad $w_8=b$\ $\checkmark$, $1$ odd
\end{minipage}
\begin{minipage}[t]{.5\textwidth}
$n=9 \rightsquigarrow (1,9)$ \quad $w_1=a$\ $\not\!\checkmark$ 

$n=9 \rightsquigarrow (3,3)$ \quad $w_3=b$\  $\checkmark$,  $3$ odd

$n=9 \rightsquigarrow (9,1)$ \quad $w_9=a$\ $\not\!\checkmark$ 
\end{minipage}

\medskip

Therefore:
\begin{equation*}
   \begin{split}
    o_w(8)=1, \ e_w(8)=1\quad\text{and}\quad 
    o_w(9)=1, \ e_w(9)=0.
   \end{split}
\end{equation*}

%
The parity functions $e_w(n)$ and $o_w(n)$ are quite irregular and we calculate their difference:
\begin{equation*}
  D_w(n)=o_w(n)-e_w(n)\,.
\end{equation*}
The difference of the divisors parity functions preserves the degree of irregularity of $e_w(n)$ and 
$o_w(n)$, taking negative, zero and positive values. 
 For each Sturmian word, $D_w(n)$ generates a path drawn by a rule similar to the odd-even drawing 
rule used to draw the face-type fractal of $w$.
Here, instead, we replace the odd-even conditions by the negative-positive sign of $D_w(n)$, 
respectively, and just draw a red dot at the point reached on the path if $D_w(n)=0$. The aspect of 
the path is always of a random walk. Two representative examples are shown in 
Figure~\ref{FigRandomWalks}.
 
 \medskip
 In order to draw information on such a irregular function, the authors of~\cite{CZ2017} 
and~\cite{CZ2019} estimated the mollified average of $D_w(n)$, which is defined by
\[M_w(x):=\sum_{n\le x}\left(1-\frac xn\right) D_w(n)\,.\]

In the following we present a sketch of the estimation of $M_w(x)$, which shows the  
asymptotic  prevalence of the odd divisors over the even ones on any Sturmian word 
(see~\cite{CZ2017} and~\cite{CZ2019} for complete details of the proof).
It turns out that $M_w(x)$ has an asymptotic behavior whose limit depends on the density:
\begin{equation*}
\beta_w:= \lim_{n\rightarrow \infty} \frac {|\{ 1\le j\le n :  w(j) = b\}|}{n}\,,
\end{equation*}
which does exists for any Sturmian words.
As an example, the density of the Fibonacci word~\cite{Dumaine2009} is
$\beta_w=\frac{3-\sqrt{5}}{2}=0.3819660\dots$.

  The estimation of $M_w(x)$ is made by first translating the word $w$  into an analytic language 
as the coefficients of the Dirichlet series 
\begin{equation*}
F(H, w, s) := \sum_{n=1}^{\infty} \frac{H(w(n))}{n^s}\,,
\end{equation*}
where $H(a) = 0$ and $H(b) = 1$. The series $F(H, w, s)$ is absolutely convergent in the half-plane 
$Re(s) > 1$ and its coefficients are related to $\beta_w$ through the limit
\begin{equation*}
\beta_w=\lim_{n\rightarrow \infty} \frac 1n\sum_{1\le j\le n} H(w(j))\,.
\end{equation*}
Then, by applying the Perron 2nd formula~\cite{Ingham}, $M_w(x)$ can be expressed as a complex  
integral:
\begin{equation*}
  \begin{split}
       M_w(x)  = 
       \frac1{2\pi i}\int_{c-i\infty}^{c+i\infty} \frac{ \left(1-\frac1{2^{s-1}}\right) 
                      \zeta(s)  F(H, w, s) x^s}{s(s+1)}\, ds, \ \ x\ge 1,\ c>1\,.
  \end{split}
 \end{equation*}
The integral is estimated by moving the path of integration to the left of the pole at $s=1$ and 
applying the residue theorem~\cite{T1986}. 

\begin{theorem}\label{Theorem41} For any Sturmian word $w$ and any $\delta>0$, we have
\begin{equation}\label{eqM}
M_w(x) = \frac{\beta_w\log 2}{2}\,x + O_\delta\left(x^{\frac13+\delta}\right)\,.
\end{equation}
\end{theorem}
In conclusion, since the main term in the estimate~\eqref{eqM} becomes positive for $x$ large 
enough, it follows that there is a significant quantifiable bias towards the odd divisors of any 
Sturmian word.

\bigskip

\textit{Acknowledgements. } 
{\small 
The author would like to thank the children who showed sacrifice and temerity while trying to 
prove uniqueness in Problem~\hyperref[labelProblem2b]{2b}, 
to R\u azvan Diaconescu and Marian V\^aj\^aitu for sharing their puzzlement, 
to Alexandru Zaharescu for his smitten ideas and inspiring solutions, to Michel Marcus and to the 
referees for carefully reading the manuscript and for their remarks on the Waring's problem, 
from Hilbert's proof to the more recent results in~\cite{DHL2000} and~\cite{Sik2016}.

Calculations and plots were generated using the free open-source 
mathematics software system \texttt{SAGE}: 
\href{http://www.sagemath.org}{http://www.sagemath.org}}.


\end{document}